\documentclass[aip,jmp,11pt,preprint]{revtex4-1}

\usepackage{amsmath,amssymb,amsthm,amsfonts}
\usepackage{graphicx}
\usepackage{dcolumn}
\usepackage{bm}
\usepackage{enumerate}
\usepackage{hyperref}
\usepackage{braket}
\usepackage[usenames,dvipsnames]{xcolor}

\newtheorem{theorem}{Theorem}[section]

\newtheorem{corollary}[theorem]{Corollary}

\newtheorem{hypo}{Hypothesis}
\newtheorem{lemma}[theorem]{Lemma}
\newtheorem{proposition}[theorem]{Proposition}
\newtheorem{example}[theorem]{Example}
\newtheorem*{remark}{Remark}
\numberwithin{equation}{section}

\newcommand{\R}{\mathbb{R}}

\newcommand{\Li}{\mathcal{L}}
\newcommand{\si}{\mathbb{S}^1}

\newcommand{\ve}{\varepsilon}

\begin{document}

\title[Adiabatic invariants of slow-fast systems]{Higher order corrections to adiabatic
invariants of generalized slow-fast Hamiltonian systems}

\author{M. Avenda\~no-Camacho}
\email{misaelave@mat.uson.mx}
\affiliation{Departamento de Matem\'aticas, Universidad de Sonora (M\'exico)}
\author{J. A. Vallejo}%
\email{jvallejo@fc.uaslp.mx}
\altaffiliation[In sabbatical leave. Permanent address: Facultad de Ciencias, UASLP (M\'exico)]{}
\author{Yu. Vorobiev}
\email{yurimv@guaymas.uson.mx}
\affiliation{Departamento de Matem\'aticas, Universidad de Sonora (M\'exico)}
\date{\today}

\begin{abstract}
We present a coordinate-free approach for constructing approximate first integrals of generalized slow-fast
Hamiltonian systems, based on the global averaging method on parameter-dependent phase spaces with $\si -$symmetry.
Explicit global formulas for approximate second-order first integrals are derived. As examples, we analyze the case quadratic in the fast variables (in particular,
the elastic pendulum), and the charged particle in a slowly-varying magnetic field.
\end{abstract}

\maketitle

\section{Introduction}

Slow-fast Hamiltonian systems appear in the theory of adiabatic approximation\cite{Arn-63,Gar-59,Kar,LW,Ne-08,SanVer-07}.
They are represented by equations of motion of the form
\begin{align}
\dot{y}=-\frac{\partial H}{\partial x},& \quad \dot{x}=\frac{\partial H}{\partial y}, \label{SF0} \\
\dot{p}=-\varepsilon\frac{\partial H}{\partial q},&\quad \dot{q}=\varepsilon\frac{\partial H}{\partial p},\nonumber
\end{align}
where $\varepsilon$ is a small perturbation parameter, and $(y,x)\in
\mathbb{R}^{2r}$, $(p,q)\in\mathbb{R}^{2k}$ are said to be fast and slow
variables, respectively. This system is Hamiltonian relative to the function
$H=H(p,q,y,x)$ and the rescaled canonical Poisson bracket on the product space
$\mathbb{R}^{2r}$ $\times\mathbb{R}^{2k}$:
\[
\{f,g\}=\left(  \frac{\partial f}{\partial y}\frac{\partial g}{\partial x}-
\frac{\partial f}{\partial x}\frac{\partial g}{\partial y}\right)
+\varepsilon\left(  \frac{\partial f}{\partial p}\frac{\partial g}{\partial q}-
\frac{\partial f}{\partial q}\frac{\partial g}{\partial p}\right) .
\]
Notice that the corresponding symplectic form
\[
\sigma=dy\wedge dx+\frac{1}{\varepsilon}dp\wedge dq ,
\]
has a singularity at $\varepsilon=0$.

Usually, systems like \eqref{SF0} appear as a result of applying a scaling
argument to an $\varepsilon -$dependent Hamiltonian on the standard phase space
$(\mathbb{R}^{2r+2k},dY\wedge dX+dP\wedge dQ)$. Two common situations arise:
\begin{enumerate}[(1)]
\item \textit{Slowly varying Hamiltonians}. In this case the original Hamiltonian has the form
\[
H(\varepsilon^{k}P,\varepsilon^{1-k}Q,Y,X),\text{ \ }(0\leq k\leq1),
\]
and the rescaling is implemented through the equations
\[
p=\varepsilon^{k}P,\text{ \ \ }q=\varepsilon^{1-k}Q,\text{ \ }%
y=\text{\ }Y,\text{ \ \ }x=X,
\]
which lead to
\[
H(\varepsilon^{k}P,\varepsilon^{1-k}Q,Y,X)=H(p,q,y,x).
\]
\item \textit{Rapidly varying Hamiltonians}. Now the original Hamiltonian is
\[
H(P,Q,\frac{Y}{\varepsilon^{k}},\frac{X}{\varepsilon^{1-k}}),\text{
\ }(0\leq k\leq1),
\]
and the rescaling
\[
p=P,\text{ \ \ }q=Q,\text{ \ \ }y=\frac{Y}{\varepsilon^{k}}\text{\ }%
,\text{\ \ }x=\frac{X}{\varepsilon^{1-k}},
\]
leads to another one in the new variables:
\[
H(P,Q,Y,X)=\varepsilon H(p,q,y,x).
\]
\end{enumerate}
In both cases, the dynamics is described by \eqref{SF0}.

These standard adiabatic models can be generalized in the following way\cite{Ku,Vor-11,VorMis-12}.
Suppose we start with two symplectic manifolds
$(M_{0},\sigma_{0})$ and $(M_{1},\sigma_{1})$. Consider the product manifold
$M=M_{0}\times M_{1}$ endowed with canonical projections $\pi_{0}:M\rightarrow M_{0}$ and $\pi
_{1}:M\rightarrow M_{1}$, and an $\varepsilon -$dependent symplectic
form%
\begin{equation}
\sigma=\pi_{0}^{\ast}\sigma_{0}+\frac{1}{\varepsilon}\pi_{1}^{\ast}\sigma_{1},
\label{SF}%
\end{equation}
for $\varepsilon\neq0$. The Poisson bracket determined by $\sigma$ is
\begin{equation}
\{,\}=\{,\}_{0}+\varepsilon\{,\}_{1}, \label{PB}%
\end{equation}
where $\{f,g\}_{0}=\Psi_{0}(df,dg)$ and $\{f,g\}_{1}=\Psi_{1}(df,dg)$
are the Poisson brackets on $M$ defined as the canonical lifts of the
non-degenerate Poisson brackets on the symplectic factors $(M_{0},\sigma_{0})$
and $(M_{1},\sigma_{1})$,\ respectively. The Poisson bivector fields $\Psi
_{0}$ and $\Psi_{1}$ on $M$ are degenerate, as $\operatorname*{rank}\Psi_{0}=\dim M_{0}$,
and $\operatorname*{rank}\Psi_{1}=\dim M_{1}$,
and the corresponding symplectic leaves are given by the slices $M_{0}\times\{m_{1}\}$
and $\{ m_{0}\}\times\{M_{1}\}$. If $(p,q,y,x)$ denote a local coordinate system on
$M$ associated to the Darboux coordinates $(p,q)$ on \ $(M_{1},\sigma_{1})$
and $(y,x)$ on $(M_{0},\sigma_{0})$, then:
\[
\Psi_{0}=\frac{\partial}{\partial y}\wedge\frac{\partial}{\partial
x},\mbox{ and }\Psi_{1}=\frac{\partial}{\partial p}\wedge
\frac{\partial}{\partial q}.
\]
The Poisson brackets $\{,\}_{0}$ \ and $\{,\}_{1}$ will be called fast and slow,
respectively.

From the viewpoint of deformation theory, the adiabatic-type Poisson bracket
\eqref{PB} is nontrivial in the following sense. The bivector field determined by \eqref{PB} reads
$\Psi=\Psi_{0}+\varepsilon\Psi_{1}$,
so $\Psi_{1}$ can be viewed as an infinitesimal deformation of $\Psi_{0}$\cite{Vais}.
One can show that this deformation is nontrivial, that is, the Poisson
cohomology class of the 2-cocycle $\Psi_{1}$ is nonzero.

For every function $H\in C^{\infty}(M)$, denote by $X_{H}$ the Hamiltonian
vector field of $H$ relative to $\sigma$, $\mathbf{i}_{X_{H}}\sigma=-dH$.
Then, $X_{H}=X_{H}^{(0)}+\varepsilon X_{H}^{(1)}$, where $X_{H}^{(0)}$ and
$X_{H}^{(1)}$ are the Hamiltonian vector fields of $H$ with respect to
$\{,\}_{0}$ and $\{,\}_{1}$, respectively:
$X_{H}^{(0)}=\mathbf{i}_{dH}\Psi_{0}\mbox{ and }X_{H}^{(1)}=\mathbf{i}%
_{dH}\Psi_{1}$.
Also, denote by $d_{0}$ and $d_{1}$ the partial exterior derivatives on $M$ along
$M_{0}$ and $M_{1}$, respectively. It is clear that $d=d_{0}+d_{1}$ is the
exterior derivative on $M$ and $d_{0}^{2}=0=d_{1}^{2}$. Also, $d_{0}\circ d_{1}%
+d_{1}\circ d_{0}=0$. Then, the relations $X_{H}^{(0)}=\mathbf{i}%
_{d_{0}H}\Psi_{0}$, $X_{H}^{(1)}=\mathbf{i}_{d_{1}H}\Psi_{1}$ hold.

By a generalized slow-fast Hamiltonian system we mean a
Hamiltonian system of the form:
\begin{equation}\label{SFH}%
(M=M_{0}\times M_{1},\{,\}=\{,\}_{0}+\varepsilon\{,\}_{1},\text{ }H).
\end{equation}
The corresponding Hamiltonian vector field $X_{H}=X_{H}^{(0)}+\varepsilon
X_{H}^{(1)}$ gives rise to a perturbed dynamics, where the unperturbed
vector field $X_{H}^{(0)}$ and the perturbation vector field $X_{H}^{(1)}$ are
Hamiltonian with respect to \emph{different} Poisson structures.
As for any perturbed Hamiltonian system, it does make sense to search for adiabatic
invariants, which are related to approximate first integrals of \eqref{SFH} as $\varepsilon \to 0$.
We are interested in the existence of (additional to the Hamiltonian $H$) approximate
first integrals of \eqref{SFH} for $\varepsilon \ll 1$, in the case when the flow of $X^{(0)}_H$ is periodic.

A $C^{\infty}$ -function on $M$,
\[
F=F_{0}+\varepsilon F_{1}+\frac{\varepsilon^{2}}{2}F_{2}+...+\frac
{\varepsilon^{k}}{k!}F_{k}+O(\varepsilon^{k+1}),
\]
smoothly depending on the small parameter $\varepsilon\ll1$, is said to be an
\emph{approximate first integral} of order $k\geq0$ for
$X_{H}=X_{H}^{(0)}+\varepsilon X_{H}^{(1)}$ if
$$\mathcal{L}_{X_{H}}F=O(\varepsilon^{k+1}).$$

An approximate first integral $F$ of order $k$ is an adiabatic invariant of order $k$
of the system \eqref{SFH}, in the sense that
\[
\mid F\circ\operatorname{Fl}_{_{X_{H}}}^{t}(p)-F(p)\mid=O(\varepsilon^k),
\]
for a long time scale: $t\sim\frac{1}{\varepsilon}$ as $\varepsilon
\rightarrow 0$ (here $\operatorname{Fl}_{_{X_{H}}}^{t}$ denotes the flow of
the vector field $X_{H}$). For $k\geq 1$, the leading term $F_0$ of $F$ must be a common
first integral of $X_{H}^{(0)}$ and the averaged perturbation vector field, and it is usually defined
as the standard action (the classical adiabatic invariant) along the periodic trajectories of
$X_{H}^{(0)}$\enspace \cite{ArKN-88,Arn-63,Gar-59}. Thus, the problem is to find the further corrections
$F_1,F_2,\ldots$ Under appropriate hypotheses for the unperturbed (fast) dynamics, we show that
a slow-fast Hamiltonian system \eqref{SFH} admits an approximate first integral of arbitrary order, and derive global
coordinate-free formulas for the first and second order corrections $F_1$ and $F_2$.
Our approach is based on the averaging method on general phase spaces with $\si -$symmetry\cite{Vor-11,VorMis-12}
and results obtained elsewhere\cite{GoKnMa,MaMoRa-90,Mon-88}. The main idea is to apply a global $\si -$normalization to \eqref{SFH}
in two stages. The first step is related to the $\si -$averaging of the original $\varepsilon -$dependent
symplectic form (the Poisson bracket) and then, in the second step, to apply the canonical normalization to
the deformed Hamiltonian on the ``new'' phase space with $\si -$symmetry. This allows us to avoid the
traditional assumption on the existence of action-angle variables, and to work on domains where the
$\si -$action is not necessarily free and trivial.
We illustrate these results with a two cases of physical interest: Hamiltonians
quadratic in the fast variables (with the elastic
pendulum as a particular example), and the charged particle in a slowly varying magnetic field.

The paper is organized as follows. In the next section we state our hypotheses and main result.
Then, in section \ref{sec3} we show how to construct the desired approximate first integrals when the perturbed system
$X_{H}=X_{H}^{(0)}+\varepsilon X_{H}^{(1)}$ is in normal form relative to the $\mathbb{S}^{1}-$action
on $M$ induced by $X^{(0)}_H$. To achieve this, a set of
homological equations must be solved. In section \ref{sec4} we prove that generalized slow-fast
system satisfying our hypotheses also satisfy the conditions required to solve these homological equations.
Section \ref{sec5} is devoted to the proof of the main theorem, and the remaining sections analyze the examples.

\section{Higher Order Corrections to Adiabatic Invariants}

\begin{hypo}[Symmetry Hypothesis]\label{hypo1}
We will assume that the flow
of the unperturbed Hamiltonian vector
field $X_{H}^{(0)}$, $\operatorname{Fl}_{X_{H}^{(0)}}^{t}$, is periodic on $M$ with frequency function $\omega\in
C^{\infty}(M),$ $\omega>0$.
\end{hypo}
That means that $\operatorname{Fl}_{X_{H}^{(0)}%
}^{t+T(m)}(m)=\operatorname{Fl}_{X_{H}^{(0)}}^{t}(m)$ for all $t\in\mathbb{R}$
and $m\in M$. Here $T=\frac{2\pi}{\omega}$ is the period function. Then, the
flow of the vector field
$
\Upsilon:=\frac{1}{\omega}X_{H}^{(0)} \label{ING1}%
$
is $2\pi -$periodic and hence $\Upsilon$ is an infinitesimal generator of the
$\mathbb{S}^{1}-$action on $M$. Let us associate\cite{MisVor-12} to that
$\mathbb{S}^{1}-$action the following operators acting on the space
$\mathcal{T}_{l}^{s}(M)$ of all tensor fields on $M$ of type $(s,l)$:
the averaging operator $\braket{\cdot}:\mathcal{T}_{l}^{s}(M)\rightarrow
\mathcal{T}_{l}^{s}(M)$,
\begin{equation*}
\braket{A}:=\frac{1}{2\pi}\int_{0}^{2\pi}(\operatorname{Fl}_{\Upsilon}^{t})^{\ast
}A\, dt, \label{AV1}%
\end{equation*}
\noindent and the integrating operator $\mathcal{S}:\mathcal{T}_{l}^{s}(M)\rightarrow
\mathcal{T}_{l}^{s}(M)$,
\begin{equation*}
\mathcal{S}(A):=\frac{1}{2\pi}\int_{0}^{2\pi}(t-\pi)(\operatorname{Fl}%
_{\Upsilon}^{t})^{\ast}A\, dt. \label{AV2}%
\end{equation*}

A tensor field $A\in\mathcal{T}_{l}^{s}(M)$ is said to be invariant with
respect to the $\mathbb{S}^{1}-$action if $(\operatorname{Fl}_{\Upsilon}%
^{t})^{\ast}A=A$ \ (for all $t\in\mathbb{R}$) or, equivalently, $\mathcal{L}%
_{\Upsilon}A=0$. \ In terms of the $\mathbb{S}^{1}-$average of $A$ the
$\mathbb{S}^{1}-$invariance condition reads $\ A=\braket{A}$. In particular,
\begin{equation}\label{invariante}
\mathcal{L}_{\Upsilon}\braket{A}=0 .
\end{equation}
The integrating
operator gives solutions to the following homological equation,%
\begin{equation}
\mathcal{L}_{\Upsilon}\circ\mathcal{S}(A)=A-\braket{A}. \label{Pr1}%
\end{equation}

It is clear that the frequency function $\omega$ and the Hamiltonian $H$ are
$\mathbb{S}^{1}-$invariant. Moreover, by the period-energy relation\cite{Gor-69,BaSn-92}
for periodic Hamiltonian flows, we have the equality%
\begin{equation}
d_{0}H\wedge d_{0}\omega=0. \label{PE}%
\end{equation}
This implies that the $\mathbb{S}^{1}-$action is canonical with respect to the
fast bracket $\{,\}_{0}$, $\mathcal{L}_{\Upsilon}\Psi_{0}=0$.
On the other hand, the $\mathbb{S}^{1}-$action does not preserve the slow
Poisson bracket $\{,\}_{1}$, in general.

Property \eqref{PE} says that the $1-$form $\frac{1}{\omega}d_{0}H$ is $d_{0}-$closed.
Our next assumption strengthens this condition.
\begin{hypo}[Parameterized Momentum Map]\label{hypo2}
We assume that
\begin{equation}
\frac{1}{\omega}d_{0}H\text{ \ \textit{is} \ }d_{0}\text{-\textit{exact},}
\label{MM1}%
\end{equation}
that is, there exists a smooth function $J:M\rightarrow\mathbb{R}$ such
that
\begin{equation}
\frac{1}{\omega}d_{0}H=d_{0}J \label{GR1}%
\end{equation}
and hence
\begin{equation}
\Upsilon=\mathbf{i}_{dJ}\Psi_{0}\text{.} \label{MM2}%
\end{equation}
\end{hypo}
This means that the $\mathbb{S}^{1}-$action associated to the periodic flow of
$X_{H}^{(0)}$ is Hamiltonian relative to $\{,\}_{0}$, with momentum map $J$.
Also, notice that condition \eqref{MM1} holds whenever $M_{0}$ is simply connected.

The momentum map $J$ is
uniquely determined by (\ref{MM2}) up to a transformation of the form
\begin{equation}
J\mapsto J+c\circ\pi_{1}, \label{FRE}%
\end{equation}
where $c$ is an arbitrary smooth function on $M_{1}$. Pick a $J$ and consider
the $1-$form $d_{1}J$ on $M$ and its $\mathbb{S}^{1}-$average $\braket{d_{1}J}$.
\begin{lemma}
There exists a closed $1-$form $\zeta\in\Omega^{1}(M_{1})$ such that
\[
\braket{d_{1}J}=\pi_{1}^{\ast}\zeta.
\]
The De Rham cohomology class $[\zeta]$ is independent of the choice of the function $c$ in \eqref{FRE}.
\end{lemma}
\begin{hypo}[Adiabatic Condition\cite{MaMoRa-90,Mon-88}]\label{hypo3}
We will assume that
$
[\zeta]=0 \label{AC1}%
$,
and hence, there exists a momentum map $J\in C^{\infty}(M)$ in (\ref{MM2})
satisfying the condition
\begin{equation}
\braket{d_{1}J}=0. \label{AC2}%
\end{equation}
\end{hypo}
Notice that the identities \eqref{GR1}, \eqref{AC2}, imply the relation%
\begin{equation}
\frac{1}{\omega}d_{1}H-d_{1}J=\frac{1}{\omega}\braket{d_{1}H}. \label{GR2}%
\end{equation}
From here and \eqref{GR1}, we conclude that sufficient conditions for the differentials
$dH$ and $dJ$ to be linearly independent at a point $m\in M$, are
$\braket{d_1 H}_m\neq 0$ and $(d_0 H)_m\neq 0$.

It is clear that condition \eqref{AC2} holds in the case when $M_{1}$ is
simply connected.
\begin{remark}\cite{MaMoRa-90,Mon-88}\label{EX}
If the ``fast'' symplectic manifold $(M_{0},\sigma_{0})$ is exact, that is,
$\sigma_{0}=d\eta$
for some $\eta\in\Omega^{1}(M_{0})$, then Hypothesis \ref{hypo1} implies Hypotheses \ref{hypo2} and \ref{hypo3}.
In this case, a momentum map $J$ satisfying the adiabatic
condition \eqref{AC2} is given by the formula
\[
J=\frac{1}{\omega}\mathbf{i}_{X_{H}^{(0)}}\braket{\pi_{0}^{\ast}\eta}.
\]
In a domain where the $\si -$action is free, this formula gives the standard action\cite{Arn-63}
along the periodic trajectories of $X^{(0)}_H$.
\end{remark}

\begin{example}
Let $\ M_{0}=\mathbb{S}^{2}\subset\mathbb{R}^{3}=\{\mathbf{x}=(x^{1}%
,x^{2},x^{3})\}$ is the unit sphere equipped with standard area form. Consider
the Hamiltonian\ on $M=\mathbb{S}^{2}\times M_{1}$  of the form
\[
H=h+\omega\mathbf{n}\cdot\mathbf{x}%
\]
for some smooth functions $h:M_{1}\rightarrow\mathbb{R}$, $\omega :M_1 \to \mathbb{R}$ and $\mathbf{n}%
:M_{1}\rightarrow \mathbb{S}^{2}$. The corresponding fast Hamiltonian system
satisfies all hypotheses above and the associated $\mathbb{S}^{1}-$action is
given by the rotations in $\mathbb{R}_{\mathbf{x}}^{3}$ about the axis
$-\mathbf{n}$. This action admits the momentum map $J=$ $\mathbf{n}%
\cdot\mathbf{x}$, which satisfies \eqref{AC2}.
\end{example}

Now, we state our main result.
\begin{theorem}\label{mainthm}
Consider a slow-fast Hamiltonian system (\ref{SFH}) and assume that the
unperturbed (fast) Hamiltonian vector field $X_{H}^{(0)}$ satisfies the
hypotheses \ref{hypo1}--\ref{hypo3} above. Let $J\in C^{\infty}(M)$ be the momentum
map of the associated $\mathbb{S}^{1}$-action, as in \eqref{MM2}, satisfying the adiabatic
condition \eqref{AC2}. Then, for any $k\geq0$, the perturbed Hamiltonian
vector field $X_{H}=X_{H}^{(0)}+\varepsilon X_{H}^{(1)}$ admits an approximate
first integral $F\in C^{\infty}(M)$, of order $k$, which is $\varepsilon
$-close to $J$. In particular, for $k=2$, the formulae
\begin{equation}\label{AI1}%
F=J+\varepsilon F_{1}+\frac{\varepsilon^{2}}{2}F_{2},
\end{equation}
where
\begin{equation}
F_{1}:=-\frac{1}{\omega}\left( \mathcal{S}(\{H,J\}_{1})+\frac{1}{2}%
\mathbf{i}_{dH}\braket{\mathbf{i}_{\mathcal{S}(d_{1}J)}\Psi_{1}}\right)  ,
\label{AI2}%
\end{equation}%
\begin{align}
F_{2}: &= \frac{2}{\omega}\mathcal{S}\left(  \{H,\frac{1}{\omega}\mathcal{S}%
(\{H,J\}_{1}\mathcal{)}\}_{1}\right. \nonumber \\
&+ \left. \frac{1}{2}\{H,\mathbf{i}_{dH}\braket{\mathbf{i}%
_{\mathcal{S}(d_{1}J)}\Psi_{1}}\}_{1}\right)  \label{AI3}%
\end{align}
give an approximate first integral of order $2$,
\[
\mathcal{L}_{X_{H}}F=O(\varepsilon^{3}).
\]
\end{theorem}
Let us remark once again that the expressions \eqref{AI1}, \eqref{AI2}, \eqref{AI3}
are intrinsic (coordinate-free) and global. They involve just the $\si -$action
induced by $X^{(0)}_H$, the Hamiltonian $H$, the ``slow'' Poisson bracket,
and the averaging operators.

The following result follows from the remark above.
\begin{corollary}
In the exact case, the assertions of Theorem \ref{mainthm} remain true
under just the symmetry hypothesis for $X_{H}^{(0)}$.
\end{corollary}
Also, from \ref{mainthm} we get the following result\cite{Arn-63,ArKN-88}.
\begin{corollary}
If $\dim M_{0}=2=\dim M_{1}$, $\braket{d_{1}H}\neq0$ and $d_0 H\neq 0$ on $M$, then, the
$2-$dimensional generalized slow-fast Hamiltonian system satisfying hypothesis \ref{hypo1}--\ref{hypo3}
is approximately integrable up to arbitrary order in $\varepsilon$.
\end{corollary}

\section{Approximate First Integrals via Normal Forms}\label{sec3}

Suppose we are given, on a manifold $M$, a perturbed vector field of the form
$A=A_{0}+\varepsilon A_{1}$. Then, a $C^{\infty}$-function on $M$,
\[
F=F_{0}+\varepsilon F_{1}+\frac{\varepsilon^{2}}{2}F_{2}+...+\frac
{\varepsilon^{k}}{k!}F_{k},
\]
is an approximate first integral of order $\varepsilon^k$ for $A$ if and
only if the functions $F_{0},...,F_{k}$ are solutions to the following
homological equations\cite{ArKN-88,Cush-84}:
\begin{align}
\mathcal{L}_{A_{0}}F_{0}& =0,\nonumber \\
\mathcal{L}_{A_{0}}F_{1}& =-\mathcal{L}_{A_{1}}F_{0},\nonumber \\
\mathcal{L}_{A_{0}}F_{2}& =-2\mathcal{L}_{A_{1}}F_{1}, \label{homological} \\
                        &\cdots \nonumber \\
\mathcal{L}_{A_{0}}F_{k}& =-k!\mathcal{L}_{A_{1}}F_{k-1}.\nonumber
\end{align}
Assume that the flow of the unperturbed vector field $A_{0}$ is
periodic with frequency function $\omega:M\rightarrow\mathbb{R}$,
$\omega>0$, and consider the $\mathbb{S}^{1}-$action on $M$ with
infinitesimal generator $\Upsilon:=\frac{1}{\omega}A_{0}$.

\begin{proposition}
Assume that there exists a smooth function $J\in C^{\infty}(M)$ which is a common
first integral of the unperturbed vector field $A_{0}$ and the
\emph{averaged} perturbation vector field $\braket{A_{1}}$, that is,
$\mathcal{L}_{A_{0}}J=0$ and $\mathcal{L}_{\braket{A_{1}}}J=0$.
Then:
\begin{enumerate}[(a)]
\item For an arbitrary $C_{1}\in C^{\infty}(M)$, the function
$$
F=J+\varepsilon(F_{1}^{0}-\braket{C_{1}}), \label{C3}%
$$
where $F_{1}^{0}:=-\frac{1}{\omega}\mathcal{L}_{\mathcal{S}(A_{1})}J$ satisfies
$\braket{F_{1}^{0}}=0$,
is an approximate second order first integral of $A$.
\item If $C_{1}$ can be chosen so that the \emph{normalization condition} of
second order holds:
\begin{equation}
\braket{\mathcal{L}_{A_{1}}F_{1}^{0}}=\braket{\mathcal{L}_{\braket{A_{1}}}C_{1}}, \label{CON}%
\end{equation}
then, $A$ admits an approximate first integral $F$, of order $2$, of the form
$$
F=J+\varepsilon(F_{1}^{0}-\braket{C_{1}})+\frac{\varepsilon^{2}}{2}(F_{2}^{0}-\braket{C_{2}}),
$$
where
$$
F_{2}^{0}:=\frac{2}{\omega}\left(  \mathcal{S}\circ\mathcal{L}_{A_{1}}\left(
\frac{1}{\omega}\mathcal{L}_{\mathcal{S}(A_{1})}J\right)  -\mathcal{L}%
_{\mathcal{S}(A_{1})}\braket{C_{1}}\right),
$$
satisfies $\braket{F_{2}^{0}}=0$, and $C_{2}\in C^{\infty}(M)$ is an arbitrary smooth function.
\end{enumerate}
\end{proposition}
\begin{proof}
The statement follows from the solvability condition for homological equations
\eqref{homological}\cite{MisVor-12,Cush-84}.
\end{proof}

Let us remark that the direct approach for finding approximate first integrals\cite{Gar-59}
leads to the verification of normalization condition \eqref{CON}, which is not
necessarily satisfied in general. It is precisely here that our hypotheses come
into play. As we will see, they will allow us first, to construct a normal form for the
perturbed Hamiltonian vector field $X_{H}=X_{H}^{(0)}+\varepsilon
X_{H}^{(1)}$, and then, from a first integral of the truncated normal form (which, in turn, is built out of
the normalization transformation and the momentum map),
an approximate first integral of the original system.

\begin{proposition}\cite{MisVor-12}
Under the symmetry hypothesis for $A_{0}$, for any $k>1$, the perturbed
vector field $A=A_{0}+\varepsilon A_{1}$ admits an $\mathbb{S}^{1}$-invariant
global normal form of order $k$,
$$
\mathcal{T}_{\varepsilon}^{\ast}(A_{0}+\varepsilon A_{1})=A_{0}+\varepsilon
\braket{A_{1}}+\frac{\varepsilon^{2}}{2}\bar{A}_{2}+...+\frac{\varepsilon^{k}}{k!}\bar{A}_{k}+O(\varepsilon^{k+1}),
$$
where $\braket{\bar{A}_i}=\bar{A}_{i}$, for $i\in\{ 2,...,k\}$.
\end{proposition}

The near-identity transformation $\mathcal{T}_{\varepsilon}$ is given by
the time$-\varepsilon$ flow of an $\varepsilon -$dependent vector field on
$M$. Therefore, the normalization transformation $\mathcal{T}_{\varepsilon}$
is well-defined on any relatively compact open subset in $M$ for small enough
$\varepsilon$. The $\mathbb{S}^{1}-$invariant vector fields $\bar{A}%
_{2},...,\bar{A}_{k}$ are well-defined on the whole $M$, and they determine the
truncated normal form\cite{MisVor-12} of order $k$ of $A$.

This leads to the following criterion.
\begin{corollary}
Suppose that an $\mathbb{S}^{1}-$invariant smooth function $J:M\rightarrow
\mathbb{R}$ is a first integral of the truncated $\mathbb{S}^{1}-$invariant
normal form of order $k$ for $A$,
\begin{equation}\label{CR}
\mathcal{L}_{\braket{A_1}}J=\mathcal{L}_{\bar{A}_{2}}J=...=\mathcal{L}_{\bar{A}_{k}}J=0.
\end{equation}
Then, the function
\[
F=J\circ\mathcal{T}_{\varepsilon}^{-1}=
J+\varepsilon F_{1}+\frac{\varepsilon^{2}}{2}F_{2}+...+\frac{\varepsilon^{k}}{k!}F_{k}+O(\varepsilon^{k+1}),
\]
where, for $s\in\{ 1,...,k\}$,
$F_{s}=\left. \frac{d^s}{d\varepsilon^s}\right|_{\varepsilon=0}(J\circ\mathcal{T}_\varepsilon^{-1})$,
is an approximate first integral on $M$ of order $k$ for $A$.
\end{corollary}

\section{Normalization of Slow-Fast Hamiltonian Systems}\label{sec4}

In this section, we show that conditions \eqref{CR} are satisfied by generalized slow-fast
Hamiltonian systems \eqref{SFH}.
To this end, we will assume that the unperturbed Hamiltonian vector field $X_{H}^{(0)}$ satisfies
the hypotheses \ref{hypo1}, \ref{hypo2}, \ref{hypo3}.
We also suppose that a momentum map $J$ as in \eqref{MM2} is
given, and that it satisfies the adiabatic condition \eqref{AC2}. As mentioned
above, generally the symplectic form $\sigma$ \eqref{SF} and the Poisson bracket \eqref{PB} are
not invariant with respect to the $\mathbb{S}^{1}-$action on $M$ associated to
the periodic flow of $X_{H}^{(0)}$. This raises the question of the
normalization of the perturbed Hamiltonian vector field $X_{H}=X_{H}^{(0)}+\varepsilon X_{H}^{(1)}$
relative to the $\mathbb{S}^{1}-$action. To
get $\mathbb{S}^{1}-$invariant global normal forms, one can apply to $X_{H}$
a non-canonical (non-Hamiltonian) Lie transform method\cite{MisVor-12}. But our point
is to maintain the Hamiltonian setting. For that purpose,
we proceed the normalization procedure in two steps. In the first stage, we
correct the drawbacks of our $\varepsilon -$dependent phase space by constructing
an $\mathbb{S}^{1}-$invariant Poisson bracket $\{,\}^{\mathrm{inv}}$ (a
symplectic form) which is $\varepsilon -$close to the original one $\{,\}$. In
a second stage, we normalize the deformed Hamiltonian $H\circ\Phi_{\varepsilon}$
up to desired order in $\varepsilon$ by applying a near-identity
transformation which will be defined as the Hamiltonian flow relative to
$\{,\}^{\mathrm{inv}}$. Here $\Phi_{\varepsilon}$ is a Poisson
isomorphism between $\{,\}^{\mathrm{inv}}$ and $\{,\}$. Thus, this
second step is basically a modification of the Deprit algorithm\cite{Dep-69} for
$\varepsilon -$dependent phase spaces.

Let us associate to the momentum map $J$ the $1-$form (on
$M$) $\Theta:=\mathcal{S}(d_{1}J)$, which has the property $\braket{\Theta}=0$.
\begin{lemma}\label{lema1}
The average $\braket{\sigma}$ of the symplectic form \eqref{SF} has the
representation $\braket{\sigma}=\sigma-d\Theta$.
Moreover, for any $\mathbb{S}^{1}-$invariant, relatively compact $N\subset M$,
and small enough $\varepsilon$, the averaged 2-form $\braket{\sigma}$ is
non-degenerate on $N$, and there exists a near-identity transformation
$\Phi_{\varepsilon}:N\rightarrow M,$ $(\Phi_{0}=\mathrm{id})$ which is
a symplectomorphism between $\braket{\sigma}$ and $\sigma$,
$\Phi_{\varepsilon}^{\ast}\sigma=\braket{\sigma}$.
\end{lemma}
\begin{proof}
The proof of this statement is based on a parametric version of the Moser
homotopy method\cite{GLS-96,DaVo-08,Vor-11}, where the symplectomorphism $\Phi_{\varepsilon}$ is
constructed as follows\cite{DaVo-08,Vor-11}.
Fix an $\mathbb{S}^{1}-$invariant, relatively compact subset $N\subset M$.
Then, define the family of 2-forms on $M$ depending on the parameter $\lambda$
\[
\delta_{\lambda}:=d_{1}\Theta+\frac{(1-\lambda)}{2}\{\Theta\wedge\Theta\}_{0},
\]
where $\frac{1}{2}\{\Theta\wedge\Theta\}_{0}$ denotes the $2-$form defined by its action
on a pair of vector fields $X,Y\in\mathcal{X}(M)$ by
$$
\frac{1}{2}\{\Theta\wedge\Theta\}_{0}(X,Y):=\{\Theta (X),\Theta (Y)\}_0,
$$
and consider the time-dependent vector field $Z_{\lambda}$ on $N$, depending on
$\varepsilon$ as a parameter, and uniquely determined by the relations
\begin{align}
\mathbf{i}_{Z_{\lambda}}\left(  \pi_{1}^{\ast}\sigma_{1}-\varepsilon
(1-\lambda)\delta_{\lambda}\right)  &=-\Theta, \label{G1}\\
\mathbf{i}_{Z_{\lambda}}\pi_{0}^{\ast}\sigma_{0}&=0. \label{G2}
\end{align}
It follows that the flow of $\varepsilon Z_{\lambda}$ is well-defined on $N$,
for small enough $\varepsilon$ and all $\lambda\in [0,1]$. Then, it suffices to take
\begin{equation}
\Phi_{\varepsilon}=\operatorname{Fl}_{\varepsilon Z_{\lambda}}^{\lambda}%
\mid_{\lambda=1} \label{PNI}%
\end{equation}
\end{proof}

Now, denote by $\{f,g\}^{\mathrm{inv}}=\Psi^{\mathrm{inv}}(df,dg)$
the non-degenerate Poisson bracket on $N$ associated to the averaged
symplectic form $\braket{\sigma}$. Then, its Poisson bivector field
$\Psi^{\mathrm{inv}}$ is $\mathbb{S}^{1}-$invariant, and it has the
representation\cite{Vor-11,VorMis-12,VorMis-13}
\begin{equation}\label{PBE1}
\Psi^{\operatorname*{inv}}=\Psi_{0}+\varepsilon\left(  \braket{\Psi_{1}}
+\mathcal{L}_{\braket{V}}\Psi_{0}\right)  +O(\varepsilon^{2}),%
\end{equation}
where%
\begin{equation}\label{PBE2}
V:=\frac{1}{2}\mathbf{i}_{\Theta}\Psi_{1}.%
\end{equation}
The following observation\cite{Vor-11} shows the role of adiabatic condition
\eqref{AC2}.
\begin{lemma}\label{lema2}
Let $J$ \ be the momentum satisfying adiabatic condition \eqref{AC2}. Then,
\begin{equation}\label{HAM}
\Upsilon=X_{J}^{(0)}=\mathbf{i}_{dJ}\Psi^{\mathrm{inv}},%
\end{equation}
that is, the $\mathbb{S}^{1}-$action is canonical on
$(N,\{,\}^{\mathrm{inv}})$ with momentum map $J$.
\end{lemma}

\begin{remark}
The adiabatic condition \eqref{AC2} was introduced\cite{MaMoRa-90,Mon-88} in the context of the theory of Hannay-Berry connections
on fibred phase spaces with symmetry. In our case, the $1-$form $\Theta =\mathcal{S}(d_1 J)$ appearing in formulas
\eqref{AI2} and \eqref{AI3}, just represents the Hamiltonian form of the Hannay-Berry connection on the trivial bundle
$M_0\times M_1 \to M_1$ over the ``slow'' base $M_1$ with the ``fast'' fiber $M_0$ endowed with an $\si -$action
associated to the periodic Hamiltonian flow of $X^{(0)}_H$\enspace\cite{GoKnMa,Vor-11,VorMis-12}. If $\Theta$ is closed, then by Lemma
\ref{lema1} the original symplectic form \eqref{SF} is $\si -$invariant.
\end{remark}

We arrive at the following normalization result.
\begin{theorem}
Suppose that the unperturbed Hamiltonian vector field $X_{H}^{(0)}$ satisfies
the hypotheses \ref{hypo1}, \ref{hypo2}, \ref{hypo3}. Then, for any $m\geq 1$ and small enough
$\varepsilon$, there exists a near-identity transformation
$\mathcal{T}_{\varepsilon}:N\rightarrow M$ which takes the original slow-fast Hamiltonian
system \eqref{SFH} into an $\mathbb{S}^{1}-$invariant normal form of order
$m$ of the form
$(N,\{,\}^{\mathrm{inv}},\text{ }H\circ\mathcal{T}_{\varepsilon})$, where
\begin{equation}\label{FC3}%
H\circ\mathcal{T}_{\varepsilon}=H+\sum_{s=1}^m\frac{\varepsilon^{s}}{s!}\braket{K_{s}}+O(\varepsilon^{m+1}),
\end{equation}
for some smooth functions $K_{1},...,K_{m}$ on $M$. In particular,
\begin{equation}\label{FC1}%
K_{1}:=\frac{1}{2}\mathbf{i}_{dH}\mathbf{i}_{\Theta}\Psi_{1}.
\end{equation}
\end{theorem}

\begin{proof}
In the first step, after applying the near-identity transformation
$\Phi_{\varepsilon}$ \eqref{PNI} to the original system \eqref{SFH}, we get a
new perturbed Hamiltonian system relative to the $\mathbb{S}^{1}-$invariant
Poisson bracket, $(N,\{,\}^{\mathrm{inv}},H\circ\Phi_{\varepsilon})$, and a deformed Hamiltonian,
\begin{equation}\label{PM2}%
H\circ\Phi_{\varepsilon}=H+\sum_{s=1}^m\frac{\varepsilon^{s}}{s!}\tilde{H}_{s}+O(\varepsilon
^{m+1}).
\end{equation}
Here the functions
\[
\tilde{H}_{s}=\frac{d^{s}}{d\varepsilon^{s}}%
\mid_{\varepsilon=0}(H\circ\Phi_{\varepsilon}), \quad s\in\{1,\ldots ,m\},
\]
are not necessarily $\mathbb{S}^{1}-$invariant. By applying Lemmas
\ref{lema1} and \ref{lema2}, a direct computation gives the result%
\begin{equation}\label{PM3}%
\tilde{H}_{1}=\frac{1}{2}\mathbf{i}_{dH}\mathbf{i}_{\Theta}\left(
\braket{\Psi_{1}}+\mathcal{L}_{\braket{V}}\Psi_{0}+\Psi_{1}\right) ,
\end{equation}
where $K_{1}$ is given by \eqref{FC1}. In a second step, we normalize
the truncated Taylor series (the $m-$th jet) of the deformed Hamiltonian
\eqref{PM2} by using the Hamiltonian flows relative to the $\mathbb{S}^{1}-$invariant
Poisson bracket $\{,\}^{\mathrm{inv}}$. Thus, we are
looking for some smooth functions on $M$, $G_{0},...,G_{m-1}$, such that the
time$-\varepsilon$ flow of the Hamiltonian vector field
$\mathbf{i}_{dG}\Psi^{\mathrm{inv}}$ of the function
\[
G=G_{0}+\sum_{i=1}^{m-1}\frac
{\varepsilon^i}{i!}G_i,
\]
gives us the desired normalization transformation,%
$$
\left( H+\sum_{s=1}^m\frac{\varepsilon^{s}}{s!}\tilde{H}_{s}\right)\circ\mathrm{Fl}%
_{\mathbf{i}_{dG}\Psi^{\mathrm{inv}}}^{\varepsilon}
=H+\sum_{s=1}^m\frac{\varepsilon^{s}}{s!}\braket{K_{s}}+O(\varepsilon^{m+1}).
$$
By a standard Lie transform argument, we conclude that $G_{0},...,G_{m-1}$,
must satisfy a set of homological equations on $M$:
\begin{align}
\{H,G_{0}\}_{0}&=\mathcal{R}_{1}-\braket{K_{1}},\nonumber \\
\{H,G_{1}\}_{0}&=\mathcal{R}_{2}-\braket{K_{2}},\label{H1} \\
&\cdots \nonumber \\
\{H,G_{m-1}\}_{0}&=\mathcal{R}_{m}-\braket{K_{m}},\nonumber
\end{align}
where the functions $\mathcal{R}_{2},...,\mathcal{R}_{m}$ are defined by the
recursive procedure given by the modified Deprit diagram. The corrections to
the standard Deprit diagram\cite{Dep-69} come from the Taylor expansion of the
Poisson tensor $\Psi^{\mathrm{inv}}$ in $\varepsilon$ at $\varepsilon
=0$. In particular, by using \eqref{PBE1}, \eqref{PBE2}, one can show that%
\begin{equation}
\mathcal{R}_{1}=\tilde{H}_{1}, \label{PM4}%
\end{equation}
and
\begin{align}
\mathcal{R}_2&=\tilde{H}_2+\mathcal{L}_{\mathbf{i}_{dG_0}\Psi_0}^{2}H+2\mathcal{L}_{\mathbf{i}_{dG_0}\Psi_0}\tilde{H}_1\nonumber \\
&+\mathbf{i}_{dH}\mathbf{i}_{dG_{0}}\left(  \braket{\Psi_{1}}+\mathcal{L}_{\braket{V}}\Psi_{0}\right) .\label{PM5}%
\end{align}
To assure the solvability of the homological equations \eqref{H1}, we choose
the functions $K_{1},\ldots ,K_{m}$ in such a way that%
\[
\braket{\mathcal{R}_{1}}=\braket{K_{1}},\ldots ,\braket{\mathcal{R}_{m}}=\braket{K_{m}}.
\]
In particular, taking into account that $H$ and $\Psi_{0}$ are $\mathbb{S}%
^{1}-$invariant, and the property $\braket{\Theta}=0$, we deduce from \eqref{PM3},
\eqref{PM4} and \eqref{PM5}, that $K_{1}$ is just given by \eqref{FC1}, and one
can put%
\[
K_{2}=\tilde{H}_{2}+\{\mathcal{S}(\frac{\tilde{H}_{1}}{\omega}),\tilde{H}_{1}\}_{0}.
\]
The global solutions to the homological equations
\eqref{H1} are given by the formulae\cite{MisVor-12},
\[
G_{0}=\frac{1}{\omega}\mathcal{S}(\mathcal{R}_{1}),\ldots ,
G_{m-1}=\frac{1}{\omega}\mathcal{S}(\mathcal{R}_{m}).
\]
Finally, the normalization transformation in \eqref{FC3} is defined by
\[
\mathcal{T}_{\varepsilon}=\Phi_{\varepsilon}\circ\mathrm{Fl}_{\mathbf{i}_{dG}\Psi^{\mathrm{inv}}}^{\varepsilon}.
\]
\end{proof}

\begin{corollary}\label{cor4}
Let $\mathcal{T}_{\varepsilon}$ be the normalization transformation in
(\ref{FC3}) and $\tilde{H}=H\circ\mathcal{T}_{\varepsilon}$. Then, the momentum
map $J$ is an approximate first integral of order $m$ for the Hamiltonian
vector field $\mathbf{i}_{d\tilde{H}}\Psi^{\mathrm{inv}}$,%
\[
\mathcal{L}_{\mathbf{i}_{d\tilde{H}}\Psi^{\mathrm{inv}}}J=O(\varepsilon^{m+1}).
\]
\end{corollary}

\begin{proof}
By using properties \eqref{HAM} and \eqref{FC3}, we get
\begin{align*}
\mathcal{L}_{\mathbf{i}_{d\tilde{H}}\Psi^{\mathrm{inv}}}J  &
=\{\tilde{H},J\}^{\mathrm{inv}}=\{\tilde{H},J\}_{0} =-\mathcal{L}_{\Upsilon}\tilde{H}\\
& =-\mathcal{L}_{\Upsilon}H-\sum_{i=1}^{m}\frac{\varepsilon^{i}}{i!}\mathcal{L}_{\Upsilon}\braket{K_{i}}+O(\varepsilon^{m+1})\\
&  =O(\varepsilon^{m+1}),
\end{align*}
where, in the last equation, we have used the $\mathbb{S}^1 -$invariance of $H$, and the property \eqref{invariante}.
\end{proof}

\section{Proof of the main result}\label{sec5}
Now, we have all the elements required to give a proof of Theorem \ref{mainthm}. Notice that
\[
X_{H}=(\mathcal{T}_{\varepsilon}^{-1})^{\ast}(\mathbf{i}_{d\tilde{H}}%
\Psi^{\mathrm{inv}}),
\]
and hence by Corollary \ref{cor4}, the function $J\circ\mathcal{T}_{\varepsilon}^{-1}$
is an approximate first integral of order $m$ for $X_{H}$. This proves the
first assertion of Theorem \ref{mainthm}. In particular, by taking the Taylor expansion of
second order%
\begin{equation}
J\circ\mathcal{T}_{\varepsilon}^{-1}=J+\varepsilon F_{1}+\frac{\varepsilon
^{2}}{2}F_{2}+O(\varepsilon^{3}),\label{TY1}%
\end{equation}
we deduce from \eqref{homological} that the functions $F_1$ and $F_2$ must satisfy the
homological equations
\begin{equation}
\mathcal{L}_{\Upsilon}F_{1}=-\frac{1}{\omega}\{H,J\}_{1},\label{TY2}%
\end{equation}
and
\begin{equation}
\mathcal{L}_{\Upsilon}F_{2}=-\frac{1}{\omega}\{H,F_{1}\}_{1}\label{TY3}%
\end{equation}
A general solution to the first equation is given by\cite{MisVor-12}
\[
F_{1}=-\frac{1}{\omega}\mathcal{S}(\{H,J\}_{1})-\braket{C_{1}},
\]
for a function $C_{1}\in C^{\infty}(M)$ which has to satisfy the solvability
condition for the second equation
\[
\frac{1}{\omega}\braket{\{H,C_{1}\}_{1}}=-\braket{\{H,\frac{1}{\omega}\mathcal{S}%
(\{H,J\}_{1}\}_{1}}.
\]
It is difficult to find $C_{1}$ from this equation. Instead of
following this approach, we will derive an explicit formula for $F_{1}$ by using the definition of
$\mathcal{T}_{\varepsilon}$. Firstly, one can verify by a direct computation
that the second term in the right-hand side of \eqref{TY1} is given by%
\[
F_{1}=-\{G_{0},J\}_{0}-\frac{1}{2}\mathbf{i}_{dJ}\mathbf{i}_\Theta%
\Psi_{0},
\]
where $G_{0}$ is a solution to the homological equation \eqref{H1}. From this fact,
and properties \eqref{GR1}, \eqref{GR2}, we get the following result.
\begin{lemma}
The second term in the Taylor expansion \eqref{TY1} is represented as follows
\[
F_{1}=-\frac{1}{\omega}\mathcal{S}(\{H,J\}_{1})-\frac{1}{\omega}\braket{K_{1}},
\]
where $K_{1}$ is given by (\ref{FC1}).
\end{lemma}

Therefore, one can put $C_{1}=\frac{1}{\omega}K_{1}$. Finally, a
\ particular solution to (\ref{TY1}) is given by the formula $F_{2}=\frac
{1}{\omega}\mathcal{S}(\{H,F_{1}\}_{1})$, which leads to the representation
(\ref{AI3}). This ends the proof of Theorem \ref{mainthm}.

\section{Hamiltonians quadratic in the fast variables}

Let us particularize the previous developments in the case of a slow-fast Hamiltonian system
$(\R^2_{(y,x)}\times \R^{2k}_{(p,q)},\frac{1}{\varepsilon}dp \wedge dq + dy\wedge dx, H)$ where
the Hamiltonian $H$ is a quadratic  function in the ``fast'' variables $\mathbf{z}=(y,x)$. To this end,
let us associate to every matrix-valued function $\mathbf{A}\in \mathfrak{sl}(2,\R )\otimes C^\infty(\R^{2k})$ the function $Q_{\mathbf{A}}=-\frac{1}{2}\mathbf{JAz\cdot z}$,
where $\displaystyle \mathbf{J}=\begin{pmatrix}
0 & -1 \\
1 &\phantom{-}0
\end{pmatrix}$, and the dot denotes the euclidian scalar product.
The Hamiltonian vector field relative to the ``fast'' Poisson bracket $\{\ ,\ \}_0$ is given by $X^{(0)}_{Q_{\mathbf{A}}}=\mathbf{Az\cdot\frac{\partial }{\partial z}}$.
Consider a Hamiltonian of the form
$
H=h + \omega Q_{\mathbf{A}},
$
for some smooth functions $h=h(p,q)$ and $\omega = \omega(p,q) >0$. We assume that $\det \mathbf{A}=1$ on an open domain in $\R^{2k}_{(p,q)}$. This implies that $X^{(0)}_H$ has periodic flow with frequency function $\omega$, hence the infinitesimal generator of the $\si -$action
is $X^{(0)}_{Q_{\mathbf{A}}}$, and the associated $\si -$action is given by the linear flow $\mathrm{Fl}^t_\Upsilon = \cos t \mathbf{I}+ \sin t \mathbf{A}$. The corresponding momentum map is $J=Q_{\mathbf{A}}$.

It is easy to see that hypotheses \ref{hypo1}--\ref{hypo3} hold in this case. For an arbitrary $\mathbf{S}\in \mathfrak{sl}(2,\R )\otimes C^\infty(\R^{2k})$,
we have the following identities,
\begin{eqnarray*}
\langle Q_{\mathbf{S}}\rangle=& \frac{1}{2}Q_{\mathbf{S-ASA}}, \\
\mathcal{S}(Q_{\mathbf{S}})=& \frac{1}{4}Q_{[\mathbf{A},\mathbf{S}]}.
\end{eqnarray*}

The Hamiltonian vector field $X_H= X_H^{(0)}+\varepsilon X_H^{(1)}$ admits an approximate first integral $F$ of second order, $\Li _{X_H} F = O(\varepsilon^3)$, of the form $F = J + \varepsilon F_1 +\frac{\varepsilon^2 }{2} F_2 +O(\varepsilon^3)$,
where $J$ is the momentum map and $F_1$, $F_2$ are given in \eqref{AI2}, \eqref{AI3}.
In the quadratic case these formulae are reduced to:
\begin{align*}
F_1=&-\frac{1}{4\omega}(Q_{[\mathbf{A,B}]}+Q_{\mathbf{A}}Q_{[\mathbf{A,C}]})\\
&+\frac{\omega}{4}\sum_{i=1}^k \left( Q_{\mathbf{A}\frac{\partial \mathbf{A}}{\partial p^i}}Q_{\frac{\partial \mathbf{A}}{\partial q^i}}-Q_{\mathbf{A}\frac{\partial \mathbf{A}}{\partial q^i}}Q_{\frac{\partial \mathbf{A}}{\partial p^i}}\right) ,\\
&=-Q_{\frac{1}{4\omega}([\mathbf{A,B}]+Q_{\mathbf{A}}[\mathbf{A,C}])-\frac{\omega}{4}\widehat{Q}(\mathbf{A})},
\end{align*}
and,
\begin{align*}
F_2=\frac{1}{2\omega}Q_{[\mathbf{A},\{h,\frac{1}{4\omega}([\mathbf{A,B}]+Q_{\mathbf{A}}[\mathbf{A,C}])
   -\frac{\omega}{4}\widehat{Q}(\mathbf{A}) \}_1]},
\end{align*}
where we have introduced the notations $\mathbf{B}:=\{h,\mathbf{A} \}_1$, $\mathbf{C}:=\{\omega,\mathbf{A} \}_1$ (the bracket being computed
separately for each coefficient of the matrix $\mathbf{A}$),
$$
\widehat{Q}(\mathbf{A}):=\sum_{i=1}^k \left( Q_{Q_{\mathbf{A}\frac{\partial \mathbf{A}}{\partial p^i}}\frac{\partial \mathbf{A}}{\partial q^i}}-Q_{\mathbf{A}\frac{\partial \mathbf{A}}{\partial q^i}}Q_{Q_{\mathbf{A}\frac{\partial \mathbf{A}}{\partial q^i}}\frac{\partial \mathbf{A}}{\partial p^i}} \right).
$$

\begin{example}
We start from the Breitenberger-Mueller model for the Hamiltonian of the elastic pendulum\cite{BM}:
\begin{equation*}
\tilde{H}=\frac{1}{2}(p^2_x+\omega^2_p x^2) +\frac{1}{2}(p^2_y+\omega^2_s y^2+\gamma x^2y).
\end{equation*}
With the rescaling $x=\frac{\mathrm{X}}{\omega_s}$, $y=\omega_p^2\mathrm{Y}$, we get,
\begin{equation*}
\tilde{H}=\frac{1}{2}(p^2_x+(\frac{\omega_p}{\omega_s})^2 \mathrm{X}^2) +\frac{1}{2}(p^2_y+(\omega_s \omega_p^2)^2 \mathrm{y}^2+\gamma(\frac{\omega_p}{\omega_s})^2 \mathrm{X}^2\mathrm{Y}).
\end{equation*}
Next, we introduce the parameter $\varepsilon =\frac{\omega_p}{\omega_s}$,
which according to the physical meaning of the problem, can be considered small: $\varepsilon \ll 1$. In terms of this parameter the Hamiltonian reads,
\begin{equation*}
\tilde{H}=\frac{1}{2}(p^2_x+\varepsilon^2 \mathrm{X}^2) +\frac{1}{2}(p^2_y+(\omega_s \omega_p^2)^2 \mathrm{y}^2+\gamma \varepsilon^2\mathrm{X}^2\mathrm{Y}),
\end{equation*}
which can be written as $\tilde{H}=H(p_x, \varepsilon \mathrm{X}, p_y, \mathrm{Y})$.
A further rescaling $p_x=p$, $q=\varepsilon \mathrm{X} $, $y= p_y$ and $x= \mathrm{Y}$, leads to the slow-fast Hamiltonian system:
\begin{equation*}
H=\frac{1}{2}(p^2+q^2)+\frac{1}{2}(y^2 +\Omega^2 x^2 +\gamma q^2 x ),
\end{equation*}
where $\Omega=\omega_s \omega_p^2$, and the $\varepsilon -$dependent Poisson bracket $\{ , \}=\{ , \}_0+\varepsilon\{ , \}_1$ on $\R^2(p,q)\times\R^2(y,x)$
(both brackets, $\{ , \}_0$ and $\{ , \}_1$, are the canonical one on $\mathbb{R}^2$). The Hamiltonian vector field $X_H = X_H^{(0)}+\varepsilon X_H^{(1)}$ is readily computed,
\begin{eqnarray*}
X_H^{(0)}&=&-(\Omega^2x+\frac{1}{2}\gamma q^2)\frac{\partial}{\partial y}+y\frac{\partial}{\partial x},\\
X_H^{(1)}&=&-(1+\gamma x)q\frac{\partial}{\partial p}+p\frac{\partial}{\partial q}.
\end{eqnarray*}
The flow of $X_H^{(0)}$ is periodic with constant frequency function $\Omega$. The momentum map of the $\mathbb{S}^1$-action $J $ in given by
\begin{equation*}
J(p,q,y,x)=\frac{\Omega}{2}(x+\frac{\gamma}{2\Omega^2}q^2 )^2+\frac{1}{2\Omega}y^2 ,
\end{equation*}
and, finally, the approximate first integral of second order for $X_H$ has the form
\begin{align*}
F&= J+\varepsilon \frac{\gamma}{\Omega^3}pqy \\
&+\varepsilon^2\frac{\gamma}{4\Omega^3}\left(\gamma q^2(x+\frac{\gamma}{2\Omega^2}q^2)\right. \\
&\left. (x-\frac{3\gamma}{2\Omega^2}q^2)+4(q^2-p^2)(x+\frac{\gamma}{2\Omega^2}q^2)-\frac{\gamma}{\Omega^2}q^2y^2 \right).
\end{align*}
\end{example}

\section{Charged particle in a slowly varying magnetic field}

The adiabatic invariants appearing in the motion of charged particles in a slowly varying
magnetic field are well-known in plasma physics since long ago\cite{Alfven,Garrido,Braun,LW},
they are related to physical phenomena such as magnetic traps\cite{Arn-63}. Indeed, for this problem
explicit expressions for corrections of the classical first and second adiabatic invariants
have been given up to first order\cite{Gall,Kruskal,Northrop}. 
Here, we will consider a cylindrically symmetric configuration for the magnetic field
(in cylindrical coordinates)
$$
\mathbf{B}=B_r\mathbf{u}_r+B_z\mathbf{u}_z= B\mathbf{u}_r -B\frac{z}{r}\mathbf{u}_z .
$$
Notice that $\frac{\partial B_r}{\partial r}=0$, so, rather than magnetic traps,
we will study a force-free contracting plasmoid.

Let us write the equations of motion of a charged ($e=1$) particle in this field in such a
way that the slow and fast components are apparent\cite{Gar-59}. In terms of the Clebsch
potentials
\begin{equation}\label{clebsch}
\begin{cases}
\alpha (\ve r,\theta ,\ve z)=-B\ve^2 rz, \\
\beta (\ve r,\theta ,\ve z)=\theta ,
\end{cases}
\end{equation}
we have $\mathbf{B}=\frac{1}{\ve^2}\bm{\nabla}\alpha\times\bm{\nabla}\beta$; also, we can write
$\mathbf{B}= \nabla \times \mathbf{A} $, where
$\mathbf{A}= \frac{1}{\varepsilon^2}\alpha\nabla\beta$. The Hamiltonian is given by
$H=\frac{1}{2}\|\widetilde{\mathbf{p}}-\mathbf{A}\|^2$, where $\widetilde{\mathbf{p}}$ is
the kinematical momentum. Let us apply the canonical transformation with generating function
\begin{equation}
F=\frac{S}{\varepsilon}p_2+\frac{\beta}{\varepsilon}p_1+\frac{\alpha}{\varepsilon}p_3-p_1p_3,\label{genfunc}
\end{equation}
where we have introduced another variable, $S$, such that $\frac{S}{\varepsilon}$ is the
arc-length along along field lines so, in the far-field regime $r\gg 1$ (which will be the one of interest for us)
\begin{equation}\label{las}
\frac{S}{\varepsilon}\simeq r .
\end{equation}
Notice that this transformation, being canonical, does not alter the (canonical) Poisson brackets.
From (\ref{genfunc}), we get
\begin{equation}
\begin{cases}
\beta &=\varepsilon q_1 + \varepsilon q_3\\
S &=\varepsilon q_2 \\
\alpha &=\varepsilon q_3 +\varepsilon p_1
\end{cases}\ .\label{eq2}
\end{equation}
Now, $\widetilde{\mathbf{p}}=\nabla F=\frac{p_2}{\varepsilon}\nabla S +\frac{p_1}{\varepsilon}\nabla \beta+\frac{p_3}{\varepsilon}\nabla\alpha $. So $\widetilde{\mathbf{p}}-\mathbf{A}=\frac{p_2}{\varepsilon}\nabla S +\frac{p_1}{\varepsilon}\nabla \beta+\frac{p_3}{\varepsilon}\nabla\alpha-\frac{1}{\varepsilon^2}\alpha \nabla \beta $. Substituting \eqref{eq2}:
\begin{equation*}
\widetilde{\mathbf{p}}-\mathbf{A}= \frac{p_2}{\varepsilon}\nabla S+\frac{p_3}{\varepsilon}\nabla\alpha -\frac{q_3}{\varepsilon}\nabla \beta ,
\end{equation*}
Thus
$$
\|\widetilde{\mathbf{p}}-\mathbf{A} \|^2
=p_2^2+\varepsilon^2B^2p_3^2(r^2+z^2)+\frac{q_3^2}{\varepsilon^2r^2}-2\varepsilon Bzp_2p_3.
$$
Now, we expand this expression in powers of $(p_3,q_3)$; to this end, we use
$\alpha = \varepsilon p_1$, $\beta = \varepsilon q_1$, $S = \varepsilon q_2$, which,
together with \eqref{clebsch}, \eqref{las}, gives
$z=-\frac{p_1}{\varepsilon B q_2}$ and $r=q_2$.
Thus, we get a Hamiltonian explicitly displaying slow and fast variables:
\begin{align*}
&H(q_1,\varepsilon p_1, \varepsilon q_2,p_2, q_3, p_3)= \\
&\frac{1}{2}\left[ p_2^2 +p_3^2\left(B^2(\varepsilon q_2)^2+\frac{(\varepsilon p_1)^2}{(\varepsilon q_2)^2}\right)+2\frac{(\varepsilon p_1)p_2}{(\varepsilon q_2)}p_3+\frac{q_3^2}{(\varepsilon q_2)^2} \right].
\end{align*}
After the obvious rescaling $\widetilde{p}_1=\varepsilon p1$, $\widetilde{q}_2=\varepsilon q_2 $, dropping the tildes
for simplicity, and noticing that $q_1$ is a cyclic variable (so we can take $p_1 =\lambda$, a parameter) we obtain the Hamiltonian system
\begin{align*}
&H(q_2,p_2,q_3,p_3)=\\
&\frac{1}{2}\left[ p_2^2 +p_3^2\left(B^2 q_2^2+\frac{ \lambda^2}{ q_2^2}\right)+2\lambda\frac{p_2}{ q_2}p_3+\frac{q_3^2}{ q_2^2} \right].
\end{align*}
The unperturbed Hamiltonian vector field (``fast'' variables) is given by
\begin{equation*}
X^{(0)}_H = -\frac{q_3}{q_2}\frac{\partial}{\partial p_3}+\left(\frac{B^2q_2^4+\lambda^2}{q_2^2}p_3+\lambda\frac{p_2}{q_2}\right)\frac{\partial}{\partial q_3},
\end{equation*}
and it has periodic flow with frequency function
$$
\omega=\frac{\sqrt{B^2q_2^4+\lambda^2}}{q_2^2}>0.
$$
The flow is of the infinitesimal generator of the $\si$-action $\Upsilon=\frac{1}{\omega}X^{(0)}_H$ is given by
\begin{eqnarray*}
\mathrm{Fl}^t_{\Upsilon}
\begin{pmatrix}
p_2 \\
q_2 \\
p_3 \\
q_3\\
\end{pmatrix}=
\begin{pmatrix}
p_2\\
q_2\\
(p_3+\lambda\frac{p_2}{\omega^2q_2^3})\cos(t)-\frac{q_3}{\omega q_2^2}\sin(t)-\lambda\frac{p_2}{\omega^2q_2^3}\\
(\omega q_2^2p_3+\lambda\frac{p_2}{\omega q_2})\sin(t)+q_3\cos(t)
\end{pmatrix},
\end{eqnarray*}
and the momentum map reads
\begin{equation*}
J(p_2,q_2,p_3,q_3):=\frac{1}{2\omega}\left(\frac{q_3^2}{ q_2^2}+ \left(\omega q_2p_3 +\lambda\frac{p_2}{\omega q_2^2}\right)^2\right).
\end{equation*}
A straightforward computation shows that, indeed, $J=\frac{v^2_\bot}{B}$, where $v_\bot$ is the transverse momentum.
From our previous results, we see that the Hamiltonian vector field $X_H= X_H^{(0)}+\varepsilon X_H^{(1)}$ admits an approximate first integral
$F$ to any arbitrary order in $\ve$. In particular, to second order it has the form
$
F = J + \varepsilon J_1 +\frac{\varepsilon^2 }{2} J_2 +O(\varepsilon^3)
$,
where $J$ is the momentum map, and
\begin{widetext}
\begin{equation*}
J_1(p_2,q_2,p_3,q_3)=-\frac{q_3 }{q_2^{10}\omega^5}
\left( p_2p_3q_2^9B^4+ \lambda q_2^4(\lambda^2p_3^2-2p_2^2q_2^2)B^2+\lambda^3(q_3^2+(p_2q_2+\lambda p_3)^2 ) \right).
\end{equation*}
\end{widetext}
The explicit expression for $J_2$ is too long to be displayed here (it can be easily computed with the aid of
a computer algebra system such as Maxima).

\end{document}